
\input amstex

\documentstyle{amsppt}

\loadbold

\magnification=\magstep1

\pageheight{9.0truein}
\pagewidth{6.5truein}



\def\And{1}
\def\BlaMon{2}
\def\CurRei{3}
\def\EtiGelone{4}
\def\EtiGeltwo{5}
\def\Kap{6}
\def\Kasetal{7}
\def\Letone{8}
\def\Lettwo{9}
\def\Mon{10}
\def\MonWit{11}
\def\Rob{12}

\topmatter

\title Commutator Hopf subalgebras and irreducible representations \endtitle

\rightheadtext{Commutator Hopf Subalgebras}

\author Edward S. Letzter \endauthor

\abstract S. Montgomery and S. Witherspoon proved that upper and lower
semisolvable, semisimple, finite dimensional Hopf algebras are of Frobenius
type when their dimensions are not divisible by the characteristic of the base
field. In this note we show that a finite dimensional, semisimple, lower
solvable Hopf algebra is always of Frobenius type, in arbitrary
characteristic.
\endabstract

\date letzter\@temple.edu \enddate

\address Department of Mathematics, Temple University, Philadelphia, PA 19122
\endaddress

\thanks This research was supported in part by a grant from the National
Security Agency. \endthanks

\endtopmatter

\document

\baselineskip = 13pt plus 2pt

\lineskip = 2pt minus 1pt

\lineskiplimit = 2pt

\head 1. Introduction \endhead 

\subhead 1.1 \endsubhead Kaplansky conjectured in \cite{\Kap} that a finite
dimensional semisimple Hopf algebra $H$ over an algebraically closed field $k$
must be of {\sl Frobenius type\/} (i.e., the dimensions of the irreducible
representations of $H$ must divide $\dim H$). Montgomery and Witherspoon
established Kaplansky's conjecture when $H$ is upper or lower semisolvable
(cf\. \cite{\MonWit, \S 3}), provided $\dim H$ is not divisible by the
characteristic of $k$ \cite{\MonWit}. Etingof and Gelaki proved that $H$ is of
Frobenius type when it is quasitriangular, provided either $k$ has
characteristic zero \cite{\EtiGeltwo} or $H$ is cosemisimple
\cite{\EtiGelone}.

\subhead 1.2 \endsubhead Continue, throughout this note, to let $k$ be an
algebraically closed field of arbitrary characteristic. Consider a series of
Hopf algebras
$$k = H_0 \subset H_1 \subset \cdots \subset H_{t-1} \subset H_t = H,$$
such that $H_{i-1}$ is normal in $H_i$, and such that $H_i/H_iH_{i-1}^+$ is
commutative, for all $1 \leq i \leq t$. (Normality is defined, e.g., in
\cite{\Mon, 3.4.1}, and $H_i^+$ denotes the kernel of the augmentation map.)
In \cite{\MonWit, \S 3}, $H$ is said to be {\sl lower solvable\/}. Lower
solvable Hopf algebras are obvious generalizations of group algebras of
solvable groups, and we can view each $H_{i-1}$ as a ``commutator Hopf
subalgebra'' of $H_i$.

In this note we prove:

\proclaim{1.3 Theorem} Let $H$ be a finite dimensional, semisimple,
lower solvable $k$-Hopf algebra. Then $H$ is of Frobenius
type. \endproclaim

\subhead 1.4 \endsubhead When the dimension of $H$ is not divisible by the
characteristic of $k$, the preceding result follows from the above cited work
of Montgomery and Witherspoon \cite{\MonWit, 3.4}. Of course, when $H$ is
quasitriangular, and when either $k$ has characteristic zero or $H$ is
cosemisimple, the result follows from the above cited work of Etingof and
Gelaki \cite{\EtiGelone; \EtiGeltwo}. We will see in \S 4 that Theorem 1.3
applies to cases not already covered in this earlier work.

\subhead 1.5 \endsubhead Our proof of Theorem 1.3 is based on the approach
found in \cite{\Letone; \Lettwo}.

\subhead 1.6 Notation \endsubhead (i) Vector spaces, tensor products,
algebras, Hopf algebras, representations, and homomorphisms are over $k$
unless otherwise noted.

(ii) We use $\Delta$ for coproducts, $\epsilon$ for counits, and $S$ for
antipodes. Given an element $h$ in a Hopf algebra $H$, we set $\Delta h = \sum
h_1 \otimes h_2$. Given a $k$-algebra homomorphism $\chi\colon H \rightarrow
k$, the $1$-dimensional module defined by $h.v = \chi(h)v$ for $h \in H$ and
$v \in k$ will be denoted $k_\chi$.

\subhead Acknowledgement \endsubhead My thanks to R. Guralnick and D. Passman
for helpful explanatory comments regarding (2.7). My thanks to the referee for
several helpful remarks and in particular for pointing me toward the examples
discussed in \S 4.

\head 2. Commutator Hopf Subalgebras \endhead 

This section provides the lemmas necessary to prove Theorem 1.3.

\subhead 2.1 Assumptions \endsubhead In this section we let $H$ and $K$ be
(not necessarily finite dimensional) $k$-Hopf algebras satisfying the
following properties: (1) $K$ is a normal Hopf subalgebra of $H$, (2) every
finite dimensional left or right simple $H$-module is semisimple as a
$K$-module, (3) every finite dimensional irreducible representation of
$H/HK^+$ is $1$-dimensional. (Recall that $HK^+ = K^+H$ is a Hopf ideal of
$H$, because $K$ is a normal Hopf subalgebra of $H$; see \cite{\Mon, 3.4.2}.)

\subhead 2.2 Example \endsubhead Let $G$ be a (not necessarily finite) group
with its commutator subgroup $G'$ having finite index. Setting $H = kG$ and $K
= kG'$, we can see as follows that the hypotheses of (2.1) are satisfied: To
start, $K$ is a normal Hopf subalgebra of $H$ because $G'$ is a normal
subgroup of $G$. Next, $H/HK^+$ is commutative and finite dimensional, and so
all of its irreducible representations are $1$-dimensional. Lastly, it follows
from Clifford Theory that every left or right simple $H$-module is semisimple
as a $K$-module; see, for example, \cite{\Rob, 8.1.3}.

\subhead 2.3 \endsubhead Let $A$ be a subring of a ring $B$. Recall that a
prime ideal $P$ of $B$ {\sl lies over\/} a prime ideal $Q$ of $A$ provided $Q$
is minimal over $P\cap A$ (i.e., there are no prime ideals of $A$ that both
contain $P\cap A$ and are properly contained in $Q$). 

\proclaim{2.4 Lemma} Let $P$ be a finite codimensional primitive ideal of $H$
lying over a (necessarily finite codimensional) primitive ideal $Q$ of $K$.
Then $H/(P + HQ)$ is a nonzero $H$-$K$-bimodule (necessarily with left
annihilator $P$ and right annihilator $Q$), and $H/(P + QH)$ is a nonzero
$K$-$H$-bimodule (necessarily with right annihilator $P$ and left annihilator
$Q$). \endproclaim

\demo{Proof} It follows from assumption 2.1(2) that $K/P\cap K$ is a
semisimple $k$-algebra. Therefore, since $Q$ is minimal over $P \cap K$, we
see that $Q.X = 0$ for some nonzero ideal $X$ of $K/P\cap K$. Hence, $(P+HQ).X
= 0$ in $H/P$, and so $P + HQ \ne H$. A similar argument shows that $P + QH
\ne H$. \qed\enddemo

\subhead 2.5 \endsubhead (i) Let $\chi \colon H \rightarrow k$ be a
$k$-algebra homomorphism. We will use $\theta _\chi$ to denote the
$k$-algebra automorphism of $H$ obtained by setting
$$\theta_\chi (h) = \sum \chi (h_1)h_2$$
for $h \in H$. (Note that $\theta_\chi$ is bijective because it has inverse
$\theta_{\chi^{-1}}$, where $\chi^{-1}$ is the convolution inverse of $\chi$.)

(ii) Let $\alpha$ be a $k$-algebra automorphism of $H$, and let $V$ be a left
$H$-module. We will use $^\alpha V$ to denote the left $H$-module with action
$h*v = \alpha(h).v$, where ``$.$'' denotes the original action on
$V$. Observe, if $\chi\colon H \rightarrow k$ is a $k$-algebra homomorphism,
that
$$^{\theta_\chi}V \; \cong \; k_\chi \otimes V ,$$
where the usual tensor-product $H$-module action is employed.

The proof of the next proposition follows arguments found in \cite{\Letone;
\Lettwo}.

\proclaim{2.6 Proposition} Let $P_1$ and $P_2$ be finite codimensional
primitive ideals of $H$ both lying over a primitive ideal $Q$ of $K$. Let
$V_1$ be a simple left $H$-module with annihilator $P_1$, and let $V_2$ be a
simple left $H$-module with annihilator $P_2$. Then there exists a 
$k$-algebra homomorphism $\chi\colon H \rightarrow k$ such that:

{\rm (i)} $P_2 = \theta_\chi(P_1)$.

{\rm (ii)} $V_1 \; \cong \; k_\chi \otimes V_2$.
\endproclaim

\demo{Proof} (i) By (2.4) there exists a finite dimensional $H$-$K$-bimodule
factor $M_1$ of $H$ with left annihilator $P_1$ and right annihilator $Q$.
Similarly, there exists a finite dimensional $K$-$H$-bimodule factor $M_2$ of
$H$ with left annihilator $Q$ and right annihilator $M_2$. Note that $M_1$ is
naturally a nonzero semisimple (and so projective) right $K/Q$-module and that
$M_2$ is naturally a nonzero semisimple left $K/Q$-module. Consequently, $M
\colon= M_1 \otimes _K M_2$ is a nonzero finite dimensional $H$-$H$-bimodule
factor of $H\otimes _K H$. The only possible left annihilator of $M$ in $H$ is
$P_1$, and the only possible right annihilator of $M$ in $H$ is
$P_2$. Furthermore, the only possible left annihilator in $H$ of any nonzero
left $H$-submodule of $M$ is $P_1$, and the only possible right annihilator in
$H$ of any nonzero right $H$-submodule of $M$ is $P_2$. Let $g$ denote the
image in $M$ of $1\otimes 1 \in H\otimes _KH$.  Since $g$ generates $M$ as an
$H$-$H$-bimodule, it follows that $g \ne 0$.

Next, observe that $M$ is a left $H$-module via
$$h*m = \sum h_1.m.S(h_2),$$
for $h \in H$ and $m \in M$. Moreover, $kg$ is a $1$-dimensional
$K$-submodule under $*$, isomorphic to $k_\epsilon$. So consider
$H*g$. Observe that it is a finite dimensional $H$-module factor of $L
\colon= H\otimes_{K}k_\epsilon$. However, $(HK^+).L = 0$, and we have
assumed that every finite dimensional irreducible representation of
$H/HK^+$ is $1$-dimensional. Therefore, there exists a nonzero element
$e$ of $M$ such that $ke$ is an $H$-module under the $*$-action.

We obtain a $k$-algebra homomorphism $\chi\colon H \rightarrow k$ such
that
$$h*e = \sum h_1.e.S(h_2) =   \chi(h)e,$$
for all $h \in H$. 

Notice, for $h \in H$, that
$$\multline h.e \quad = \quad \sum h_1\epsilon (h_2).e \quad = \quad \sum
h_1.e\epsilon (h_2) \quad = \quad \sum h_1.e.S(h_2)h_3 \quad \\ = \quad \sum
(h_1*e).h_2 \quad = \quad \sum e.\chi(h_1)h_2 \quad = \quad e .\theta_\chi
(h). \endmultline$$

Since $0 = P_1.e = e.\theta_\chi(P_1)$, it follows that $\theta_\chi(P_1)
\subseteq P_2$. Since $\theta_\chi(P_1)$ is a maximal ideal of $H$, it follows
that $\theta_\chi(P_1) = P_2$.  Part (i) follows.

(ii) Set $U$ equal to the left $H$-module $^{\theta_\chi}V_2 \cong k_\chi
\otimes V_2$. By (i), $P_1.U = \theta_\chi(P_1).V_2 = P_2.V_2 = 0$. Since $U$
is a simple $H$-module, it must follow that $U \cong V_1$.  \qed\enddemo

\subhead 2.7 Example \endsubhead The following example, illustrating the
applicability of (2.6), served as the original motivation for our work in this
note: Set $H = kG$ and $K = kG'$, as in (2.2). For each $g \in G$, let
$\alpha_g$ be the $k$-algebra automorphism of $H$ mapping $h \in H$ to
$g^{-1}hg$. Since $G'$ is normal in $G$, each $\alpha_g$ restricts to an
automorphism of $K$ that we will also denote by $\alpha_g$.

Now let $U$ be an irreducible finite dimensional left $K$-module, and let $Q$
be the annihilator of $U$ in $K$. Set $M = H\otimes_KU$. Note, for each $g \in
G$, that the $K$-submodule $g\otimes U$ of $M$ is isomorphic to $^{\alpha_g}U$
and so has annihilator $\alpha_{g^{-1}}(Q)$ in $K$. Of course,
$$M \; = \; \bigoplus _{g \in C} g \otimes U,$$
where $C$ is a set of coset representatives of $G'$ in $G$. Therefore, the
annihilator in $K$ of an arbitrary $K$-module composition factor of $M$ must
equal $\alpha_{g^{-1}}(Q)$, for some choice of $g \in G$.

Next, choose an $H$-module composition factor $V$ of $M$, with annihilator $P$
in $H$. It follows from the last paragraph that $P$ lies over
$\alpha_{g^{-1}}(Q)$, for some $g \in G$, and hence $\alpha_g(P)$ lies over
$Q$. Since $V$ was chosen arbitrarily, it now follows from (2.6) that there
exist automorphisms $\beta_1,\ldots,\beta_t$ of $H$ for which
$$^{\beta_1}V,\ldots,\hbox{}^{\beta_t}V$$
is a complete list of the $H$-module composition factors of $M$. In particular,
$\dim V$ divides $\dim M = [G:G']\dim U$. Of course, each finite dimensional
simple $H$-module is a composition factor of some $H\otimes_KU$, for an
appropriate choice of $U$.

The preceding conclusions can all be derived as direct consequences of
standard Clifford Theory. In particular, these conclusions can be proved
using, for instance, the arguments in \cite{\CurRei, 49.2, 51.7}. (Thanks to
R. Guralnick and D. Passman for their helpful comments.)

\proclaim{2.8 Lemma} Assume that $H$ is finitely generated as a right
$K$-module and that every finite dimensional left $H$-module is semisimple.
Let $Q$ be a finite codimensional primitive ideal of $K$, and suppose that $P$
is the annihilator in $H$ of some composition factor of the
finite dimensional left $H$-module $H \otimes _K (K/Q)$. Then $P$ lies over
$Q$. \endproclaim

\demo{Proof} Set $M = H \otimes _K (K/Q)$, and let $e$ denote the image of
$1\otimes 1$ in $M$. Since $M$ is semisimple as a left $H$-module, there
exists a simple left $H$-module factor (i.e., image) $N$ of $M$ with left
annihiltor $P$. Consequently, $M/PM$ is a nonzero $H$-$K$-bimodule factor of
$M$, and $M/PM$ must have left annihilator $P$ and right annihilator $Q$. But
$Q.e = e.Q = 0$, and so $Q$ is the annihilator of a left $K$-submodule of
$M/PM$. Hence, $Q$ is minimal over $P \cap K$. \qed\enddemo

\proclaim{2.9 Proposition} Assume that $H$ is finitely generated as a right
$K$-module, that every finite dimensional left $H$-module is semisimple, and
that $U$ is a finite dimensional simple left $K$-module.

{\rm (i)} Let $V_1,\ldots,V_s$ be the $H$-module composition factors of
$H\otimes _KU$. Then there exist $k$-algebra homomorphisms
$\chi_2,\ldots,\chi_s\colon H \rightarrow k$ such that
$$V_1 \; \cong \; k_{\chi_2}\otimes V_2 \; \cong \; \cdots \; \cong \;
k_{\chi_s} \otimes V_s.$$
Consequently, $\dim V_1 = \dim V_2 = \cdots = \dim V_s$.

{\rm (ii)} Suppose that $H$ is a free right $K$-module of rank $m$
and that $V$ is an $H$-module composition factor of $H\otimes _K
U$. Then $\dim V$ divides $m\dim U$.
\endproclaim

\demo{Proof} (i) Let $P_1,\ldots,P_s$ be the respective annihilators in $H$ of
$V_1,\ldots,V_s$. Then each $P_i$ is the annihilator of some composition
factor of $H\otimes _K(K/Q)$, where $Q$ is the annihilator in $K$ of $U$. By
(2.8), $P_1,\ldots,P_s$ all lie over $Q$. Part (i) now follows from (2.6).

(ii) Follows easily from (i). \qed\enddemo

\head 3. Proof of Theorem \endhead

This section is a proof of Theorem 1.3\ Assume throughout that $H$ is a finite
dimensional semisimple $k$-Hopf algebra and that there exists a series of Hopf
subalgebras
$$k = H_0 \subset H_1 \subset \cdots H_{t-1} \subset H_t = H,$$
such that $H_{i-1}$ is normal in $H_i$, and such that
$H_i/H_iH_{i-1}^+$ is commutative for all $1 \leq i \leq t$.

It follows from the Nichols-Zoeller theorem
(see, e.g., \cite{\Mon, 3.1.5}) that $H_j$ is free as a left and right
$H_i$-module, for all $1 \leq i \leq j \leq t$. It follows from
\cite{\Mon, 3.2.3} that $H_0, H_1, \ldots, H_t$ are all semisimple.

Now set $K = H_{t-1}$, and assume by induction that $K$ is of
Frobenius type. Let $m$ denote the rank of $H$ as a free right
$K$-module, and let $n = \dim K$. Note that $H$ and $K$ satisfy all of
the hypotheses in all of the lemmas in \S 2.

Next, let $V$ be a simple left $H$-module. There exists a simple left
$K$-module $U$ such that $V$ is an $H$-module composition factor of
$H\otimes_KU$. It follows from (2.9ii) that $\dim V$ divides $m\dim U$.  By
the induction hypothesis $\dim U$ divides $n$, and so $\dim V$ divides $mn =
\dim H$. The theorem follows. \qed

\head 4. A Class of Applicable Examples \endhead

In this section we briefly outline a class of examples covered by Theorem 1.3
but not by the earlier results \cite{\EtiGelone; \EtiGeltwo; \MonWit}
described in (1.1). Our approach is based on \cite{\And} (cf\. references
cited therein).

\subhead 4.1 \endsubhead Following \cite{\And, \S 3}, let $A$ and $B$ be
finite dimensional $k$-Hopf algebras, let $\, \rightharpoonup \, \colon
B\otimes A \rightarrow A$ be a weak action (in the sense of \cite{\And, 2.1}),
and let $\rho \colon B \rightarrow B\otimes A$ be a weak coaction (in the
sense of \cite{\And, 2.3}). For $\sigma \colon B\otimes B \rightarrow A$ and
$\tau \colon B \rightarrow A \otimes A$ satisfying certain restrictions
(detailed in \cite{\And, \S 3.1}), one can construct a Hopf algebra
$$H = A^\tau \# _\sigma B;$$
see \cite{\And, 3.1.12}. Retain the preceding notation for the remainder of
this section.

\subhead 4.2 \endsubhead It is further proved in \cite{\And, \S 3} that: (i) 
$H$ is isomorphic as an algebra to the crossed product $A\#_\sigma B$. (ii)
$A$ can be identified isomorphically with a normal Hopf subalgebra of
$H$, and
$$B \; \cong \; H/HA^+.$$

\subhead 4.3 \endsubhead Further assume, for the remainder, that $k$ has
positive characteristic $p$, that $X$ is a finite solvable group whose order
is not divisible by $p$, that $Y$ is a finite group whose order is divisible
by $p$, that $A = kX$, and that $B = (kY)^*$. (It follows that $H^*$ is an
{\sl abelian extension\/} in the sense of \cite{\Kasetal}.)

(i) Since $A$ and $B$ are semisimple, it follows from (4.2i) and
\cite{\BlaMon} that $H$ is semisimple. Since $B$ is commutative and $X$ is
solvable, it follows from (4.2ii) that $H$ is lower solvable. Consequently,
Theorem 1.3 applies.

(ii) Since the characteristic of $k$ divides $\dim H$, the approach in
\cite{\EtiGeltwo} and \cite{\MonWit} cannot be applied. (In particular, since
the characteristic of $k$ divides $\dim B$, \cite{\MonWit, 2.1} does not
apply.)  Also, $H$ is not cosemisimple because $B$ is not cosemisimple (see,
e.g., \cite{\Mon, 3.2.3}), and so the approach in \cite{\EtiGelone} cannot be
applied.

\subhead 4.4 \endsubhead More concretely, suppose that $Y$ is a finite group
whose order is divisible by $p$, and suppose that $X$ is a solvable subgroup of
$Y$ whose order is not divisible by $p$. Then $A$ is a left $B$-module algebra
via the action
$$b \, \rightharpoonup \, a \; = \; \sum a_1 \langle b, a_2 \rangle $$
for $a \in A$ and $b \in B$. Letting $H$ now denote the smash product $A\# B$,
we can see as above that Theorem 1.3 applies to $H$ but that \cite{\EtiGelone,
\EtiGeltwo, \MonWit} do not.

\enddocument